# Nonlinear renewal theorems for random walks with perturbations of intermediate order

**Keiji Nagai[1] and Cun-Hui Zhang[2],\***

*Yokohama National University and Rutgers University*

**Abstract:** We develop nonlinear renewal theorems for a perturbed random walk without assuming stochastic boundedness of centered perturbation terms. A second order expansion of the expected stopping time is obtained via the uniform integrability of the difference between certain linear and nonlinear stopping rules. An intermediate renewal theorem is obtained which provides expansions between the nonlinear versions of the elementary and regular renewal theorems. The expected sample size of a two-sample rank sequential probability ratio test is considered as the motivating example.

## 1. Introduction

Let $X, X_n, n \geq 1$, be i.i.d. random variables with a finite, positive mean $EX = \mu$. Let $\{Z_n, n \geq 1\}$ be a perturbed random walk defined as

$$(1.1) \qquad Z_n = S_n + \xi_n, \quad S_n = X_1 + \cdots + X_n, \qquad n \geq 1,$$

where $\{\xi_n, n \geq 1\}$ are random variables such that $\{X_1, \xi_1, \ldots, X_n, \xi_n\}$ are independent of $\{X_{n+j}, j \geq 1\}$ for all $n \geq 1$. Define

$$(1.2) \qquad T_b = \inf\{n \geq 1 : Z_n > b\}, \quad \tau_b = \inf\{n \geq 1 : S_n > b\}.$$

Nonlinear renewal theory studies probabilistic quantities related to stopping rules $T_b$, especially Blackwell-type theorems for the convergence of the renewal measure $U_b((b, b+h]) = \sum_n P\{b < Z_n \leq b + h\}$, the distribution of the excess over the boundary $R_b = Z_{T_b} - b$, and expansions of $ET_b$ and $\text{Var}(T_b)$, as extensions of the (linear) renewal theorems concerning probabilistic quantities related to $\tau_b$. Many important applications of nonlinear renewal theory come from sequential analysis in which nonlinear renewal theory provides crucial analytical tools and methodologies [8, 19, 22, 25, 28].

Many authors have studied nonlinear renewal theory. See for example Chow and Robbins [6], Chow [3], Siegmund [23, 24], Gut [9], Pollak and Siegmund [18], Woodroofe [26, 27], Lai and Siegmund [15, 16], Chow, Hsiung and Lai [4], Lalley [17], Hagwood and Woodroofe [10], Woodroofe and Keener [31], Zhang [32], Hu [11],

---

[1]Yokohama National Univeristy, Faculty of Economics, Yokohama 240-8501, Japan, e-mail: knagai@ynu.ac.jp
[2]Department of Statistics, Hill Center, Busch Campus, Rutgers University, Piscataway, NJ 08854, USA, e-mail: czhang@stat.rutgers.edu
*Research partially supported by National Science Foundation.

Aras and Woodroofe [1], Kim and Woodroofe [12, 13], and the books by Woodroofe [28] and Siegmund [25] on the subject. A main condition of the existing Blackwell-type nonlinear renewal theorems is

$$(1.3) \qquad \lim_{\theta \to 0} \limsup_{n \to \infty} P\left\{ \max_{1 \le j \le \theta n^\alpha} \left| \xi_{n+j} - \xi_n \right| > \epsilon \right\} = 0, \quad \forall \epsilon > 0,$$

(Lai and Siegmund [15]). This condition of slowly changing perturbation allows $\xi_n$ with unbounded variability [i.e. stochastically unbounded $\xi_n - c(n)$ for any centering constants $c(n)$]. However, in addition to (1.3) and other regularity conditions, existing results on the second order [i.e. up to $o(1)$] expansion of $ET_b$ requires

$$(1.4) \qquad \text{the uniform integrability of} \max_{N_0 \vee n \le j \le n+n^\alpha} \left| \xi_j - f(j) \right|$$

for certain $1/2 < \alpha \le 1$, where $f(t)$ is a slowly changing deterministic function and $N_0$ is a random variable with $EN_0 < \infty$ [15]. This condition precludes perturbation processes $\xi_n$ with unbounded variability. In this paper, we remove the restriction on the bounded variability of $\xi_n$ by imposing the condition of

$$(1.5) \qquad \text{the uniform integrability of} \max_{1 \le j \le Mn^\alpha} \left( |\xi_{n+j} - \zeta_n| \wedge n^\alpha \right)$$

for all $M \in \mathbb{R}$, instead of (1.4), where $\zeta_n$ are certain truncated $\xi_n$. This will be done in Section 2.

Blackwell-type nonlinear renewal theorems and second order expansions of $ET_b$ were first obtained by Woodroofe [26] for stopping rules of the form

$$(1.6) \qquad T_b = \inf\left\{ n \ge 1 : S_n > A(n;b) \right\}$$

with certain nonlinear boundaries $A(\cdot;b)$. Lai and Siegmund [15, 16] pointed out that in many applications (1.6) can be written as (1.2) for certain (possibly different) random walk with perturbation and developed nonlinear renewal theorems for (1.2) under much weaker conditions on the distribution of $X$. Zhang [32] studied nonlinear renewal theorems for both (1.2) and (1.6) through the uniform integrability of $|T_b - \tau_a|^p$, where $a = b - f(b/\mu)$ with the function $f(t)$ in (1.4). He also obtained a second order expansion of $\text{Var}(T_b)$ for stopping rules of the form (1.6), and thus demonstrated certain advantages of investigating (1.6). The methods in Zhang [32] and here can be combined to study a general form of stopping rules with (1.2) and (1.6) as special cases, but for simplicity we confine the rest of our discussion to stopping rules of form (1.2).

The motivating example of this paper is Savage and Sethuraman's [21] two-sample rank sequential probability ratio test (SPRT) for $G = F$ against a Lehmann alternative. In Section 3 we provide formulas for the expansion of the expected sample size of the rank SPRT, with an outline of a proof. Our calculation shows that the rank log-likelihood ratio is of the form (1.1) such that $(\xi_n - E\xi_n)/(\log n)^{1/2}$ converges to a nondegenerate normal distribution under the null hypothesis, so that we are truly dealing with perturbations of intermediate order.

Our expansion of $ET_b$ is obtained via the uniform integrability of $|T_b - \tau_b^*|$ under (1.5), where $\tau_b^*$, defined in (4.1) below, is a linear stopping rule with perturbation at an initial time $n = n_*$. In fact, we will develop sufficient conditions for the uniform integrability of $\{|T_b - \tau_b^*|/\rho(b)\}^p$ for certain normalization constants $1 \le \rho(b) = o(b)$ and consequently an intermediate renewal theory for the expansion of $ET_b$ up to $O(\rho(b))$. These will be done in Section 4.

## 2. Expectation of stopping rules

In this section, we provide a second order expansion of $ET_b$ under (1.5). Let $\rho(x)$ be a function satisfying

(2.1) $\quad 1 \leq \rho(x) = o(x)$ as $x \to \infty$, $\quad \sup_{x>1} \sup_{x \leq t \leq 3x} \rho(t)/\rho(2x) < \infty.$

We shall first state a general set of regularity conditions on the perturbation $\xi_n$ in (1.1) which we call $\rho$-regularity. We denote throughout the sequel $\lfloor x \rfloor$ the integer part of $x$ and $\lceil x \rceil$ the smallest integer upper bound of $x$.

The process $\{\xi_n\}$ in (1.1) is called *$\rho$-regular*, or regular of order $\rho(\cdot)$, with parameters $p \geq 1$ and $1/2 < \alpha \leq 1$ if the following three conditions hold: there exist constants $\delta_0 > 0$, $\theta > 0$, $\theta < \mu$ for $\alpha = 1$, $K > 0$, $w_0 > 0$ and $0 < \theta^* < K\mu$ such that

(2.2) $$P\Big\{\max_{\delta_0 n < j \leq n} \xi_j > \theta n^\alpha\Big\} = o(\rho(n)/n)^p,$$

(2.3) $$\sum_{k \geq n + Kn^\alpha} k^{p-1} P\{\xi_k \leq -(k-n)\mu + w_0 k^\alpha\} = o(\rho^p(n)),$$

and that for $\zeta_n = (\xi_n \wedge \theta n^\alpha) \vee (-\theta^* n^\alpha)$ and all $M < \infty$

(2.4) $$\left\{\max_{0 \leq j \leq M n^\alpha} \left(\frac{|\zeta_n - \xi_{n+j}| \wedge n^\alpha}{\rho(n)}\right)^p, n \geq 1\right\} \text{ is uniformly integrable.}$$

We impose the 1-regularity with $\rho(x) = 1$ in Theorem 2.1 below for the expansion of $ET_b$ and the general $\rho$-regularity in Section 4 for uniform integrability and an intermediate renewal theorem of order $\rho(\cdot)$.

**Theorem 2.1.** *Suppose $\{\xi_n, n \geq 1\}$ is 1-regular, i.e. $\rho(x) = 1$, with parameters $p = 1$ and $1/2 < \alpha \leq 1$. Suppose $X$ is non-lattice with $EX = \mu > 0$ and $E|X|^{2/\alpha} < \infty$. Suppose $bP\{T_b \leq \delta_0 b/\mu\} = o(1)$ for the $\delta_0$ in (2.2) and the slowly changing condition*

(2.5) $$\max_{1 \leq j \leq n^\alpha} |\xi_{n+j} - \zeta_n| = o_P(1)$$

*holds with the $\zeta_n = (\xi_n \wedge \theta n^\alpha) \vee (-\theta^* n^\alpha)$ in (2.4). Then, as $b \to \infty$,*

(2.6) $$\mu ET_b = b - E\zeta_{n_b} + \frac{ES_{\tau_0}^2}{2ES_{\tau_0}} + o(1),$$

*where $n_b = \lfloor b/\mu \rfloor$, $S_n$ is as in (1.1) and $\tau_0$ is as in (1.2).*

**Remark 2.1.** As mentioned earlier, the main difference between the $\rho$-regularity and the usual regularity conditions in the nonlinear renewal theory literature is (2.4), which allows $\xi_n$ to have unbounded variability as $n \to \infty$.

**Remark 2.2.** Since $(k-n)/k^\alpha$ is increasing in $k$, $(k-n)\mu - w_0 k^\alpha \geq (\mu K/(1+K)^\alpha - w_0)k^\alpha$ for $k \geq n + K(n\mu)^\alpha$. Thus, (2.3) is a consequence of

(2.7) $$\sum_{k=1}^\infty k^{p-1} P\{-\xi_k \geq w_1 k^\alpha\} < \infty$$

for certain $w_1 > 0$. In typical nonlinear renewal theorems, $\rho(x) = 1$ and (2.7) is imposed (with $w_1 < \mu$ for $\alpha = 1$) instead of (2.3).

**Remark 2.3.** The condition $bP\{T_b \leq \delta_0 b/\mu\} \to 0$ holds if $EX^2 < \infty$ and $nP\{\max_{j \leq n} \xi_j > w_2 n\} = o(1)$ for some $w_2 < \mu(1/\delta_0 - 1)$.

*Proof of Theorem* 2.1. Let $T = T_b$ be as in (1.2), $n_* = \lfloor b/\mu - \eta_* b^\alpha \rfloor$, and $\tau^* = \tau_b^* = \inf\{n \geq n_* : S_n + \zeta_{n_*} > b\}$ be as in (4.1) below with $\eta_* > 0$, $\eta_* < 1/\mu$ for $\alpha = 1$. Let $R_b^* = S_{\tau^*} + \zeta_{n_*} - b$ be the overshoot for $\tau^*$. Since $\zeta_{n_*}$ is bounded and $b - \zeta_{n_*} \to \infty$,

$$(2.8) \qquad \mu E\tau^* = ES_{\tau^*} = ER_b^* + b - E\zeta_{n_*} = b - E\zeta_{n_*} + \frac{ES_{\tau_0}^2}{2ES_{\tau_0}} + o(1)$$

by the Wald identity and the standard linear renewal theorem [7]. The proof is based on (2.8), (2.4) with $\rho(\cdot) = 1$, (2.5), and the uniform integrability of $|T - \tau^*|$ in Theorem 4.1.

Since $E|T - \tau^*|$ is bounded, $ET < \infty$ and $\mu E(T - \tau^*) = E(S_T - S_{\tau^*})$ by the Wald identity. Let $R_b = Z_T - b$ be the overshoot for $T$ in (1.2). By (1.2) and (4.1)

$$(2.9) \qquad S_T - S_{\tau^*} = (R_b + b - \xi_T) - (R_\tau^* + b - \zeta_{n_*}) = R_b - R_b^* + (\zeta_{n_*} - \xi_T).$$

Since $EX^2 < \infty$, the uniform integrability of $|T - \tau^*|$ implies that of $|S_T - S_{\tau^*}|$. Since $\zeta_{n_*} = O(b^\alpha)$, $|\tau^* - b/\mu| = O_P(b^\alpha)$. This and $T - \tau^* = O_P(1)$ imply

$$(2.10) \qquad |T - b/\mu| + |\tau^* - b/\mu| = O_P(b^\alpha),$$

so that $\zeta_{n_*} - \xi_T = o_P(1)$ by (2.5). Moreover, (2.10) and (2.5) imply the convergence of both $R_b$ and $R_b^*$ in univariate distribution to the same limit [15]. These facts and the uniform integrability of $|S_T - S_{\tau^*}|$ imply $\mu E(T - \tau^*) = E(S_T - S_{\tau^*}) \to 0$ in view of (2.9). It then follows from (2.8) that $ET = b - E\zeta_{n_*} + ES_{\tau_0}^2/(2ES_{\tau_0}) + o(1)$. Since $E|\zeta_{n_*} - \zeta_{n_b}| \to 0$ by (2.4) and (2.5), (2.6) follows. □

## 3. Rank SPRT

Given a constant $\Delta > 0$ and two independent samples of equal size from continuous distribution functions $F$ and $G$ respectively, the rank likelihood ratio for testing $H_0 : G = F$ against the Lehmann alternative $G = F^\Delta$ is $\Lambda_n = L_n(F, F^\Delta)/L_n(F, F)$, where $L_n(F, G)$ is the probability mass function of the ranks of the $F$-sample within the combination of the $F$- and $G$-samples. Suppose pairs of observations, one from $F$ and one from $G$, are taken sequentially, the rank SPRT [21] rejects $H_0$ iff $\Lambda_T > e^b$ with the stopping rule

$$(3.1) \qquad T = T_{a,b} = \inf\{n \geq 1 : Z_n < -a \text{ or } Z_n > b\}, \quad Z_n = \log \Lambda_n.$$

In this section, we provide formulas for the expansion of the expectation of the sample size (3.1) of the rank SPRT under the following conditions:

$$(3.2) \qquad G = F^A, \quad \mu = \log \Delta + \int \log\left(\frac{F+G}{F+\Delta G}\right) d(F+G) \neq 0.$$

We outline a proof of the expansion via a representation of the rank log-likelihood ratio $Z_n$ in (3.1) as a perturbed random walk (1.1) satisfying the conditions of Theorem 2.1. We prove in (3.9) the convergence of $(\xi_n - E\xi_n)/(\log n)^{1/2}$ in distribution under the null hypothesis, so that the rank log-likelihood is truly a random walk with perturbations of intermediate order.

The rank likelihood ratio was considered earlier by Savage [20] who showed that

$$(3.3) \qquad Z_n = \log\left(\Delta^n(2n)!\right) - n\int \log(nF_n + \Delta nG_n)d(F_n + G_n),$$

where $F_n$ and $G_n$ are the ECDF's of the $F$- and $G$-samples. Chernoff and Savage provided [2] representations of certain rank statistics as (1.1) and proved their asymptotic normality, i.e. $\xi_n = o_P(n^{1/2})$. Lai identified [14] the rank log-likelihood ratio $Z_n$ with a Chernoff–Savage statistic and proved the quick convergence $n^{-\epsilon}\xi_n$ for all $\epsilon > 0$ in (1.1). The problem was further investigated by Woodroofe [28, 29] who proved (1.3) for the perturbation term of (3.3) and derived asymptotic formulas for the error probabilities of the rank SPRT. Our results are based on a finer expansion of (3.3) from a slightly different expression (3.4) below.

Let $H = F + G$, $W = F + \Delta G$, and $H_n$ and $W_n$ be their empirical versions. Since $\log((2n)!) = n\int \log(nH_n)dH_n$, (3.3) can be written as

$$(3.4) \qquad Z_n = n\log(\Delta) + n\int \log\left(\frac{H_n}{W_n}\right)dH_n.$$

Let $\psi_n = H_n/W_n$ and $\psi = H/W$. We write (3.4) as a perturbed random walk (1.1) with

$$(3.5) \qquad S_n = n\log\Delta + n\int \log(\psi)dH_n + n\int (H_n/H - W_n/W)dH$$

and the perturbation term

$$(3.6) \qquad \xi_n = \int n\Big\{\log(\psi_n/\psi)dH_n - (H_n/H - W_n/W)dH\Big\}.$$

Let $y_k$ be hypergeometric$(n,n,k)$ variables and $\eta = (1-\Delta)/(1+\Delta)$. Define

$$C(\eta) = \lim_n \sum_{k=1}^{2n} E\log\left(1 + \eta\Big(\frac{2y_k}{k} - 1\Big)\right) + \frac{\eta^2}{2}\log(2n).$$

The above limit exists since $Ey_k = k/2$ and $\text{Var}(y_k) = (k/4)(2n-k)/(2n-1)$. We state our expansion of the expectation of the sample size (3.1) as follows.

**Theorem 3.1.** *Set $\eta = (\Delta - 1)/(\Delta + 1)$. Let $T_{a,b}$, $S_n$, and $\tau_0$ be as in (3.1), (3.5), and (1.2) respectively. Suppose $\mu > 0$ in (3.2) and $a/\log b \to \infty$ as $b \to \infty$. Then,*

$$(3.7) \qquad \mu ET_{a,b} = b - E\xi_{n_b} + ES_{\tau_0}^2/(2ES_{\tau_0}) + o(1),$$

*where $n_b = \lfloor b/\mu \rfloor$ as in Theorem 2.1. Moreover, if $G = F^A$, then*

$$\mu ET_{a,b} = \begin{cases} b - \int h(x)d(x + x^A) + ES_{\tau_0}^2/(2ES_{\tau_0}) + o(1), & A \neq 1 \\ b - (\eta^2/2)\log(2b/\mu) + C(\eta) + ES_{\tau_0}^2/(2ES_{\tau_0}) + o(1), & A = 1, \end{cases}$$

*where $h(x) = (1-\Delta)^2 x^{1+A}/\{2(x + \Delta x^A)^2(x + x^A)\}$.*

**Remark 3.1.** Since $\Lambda_n = L_n(G^{1/\Delta}, G)/L_n(G, G)$, Theorem 3.1 also provides expansions of $ET_{a,b}$ for $\mu < 0$ as $(a, b/\log a) \to (\infty, \infty)$ by symmetry.

Here is an outline of a proof of Theorem 3.1. For $\epsilon_0 \leq \mu$, $P(Z_n \leq 0) \leq P(|Z_n - n\mu| > \epsilon_0 n) \to 0$ exponentially fast [21]. Since $P(|Z_n| \leq n|\log \Delta|) = 1$, this implies $P(Z_{T_{a,b}} < 0) = o(b^2)$, so that $ET_{a,b} = ET_b + o(1)$. Thus, it suffices to consider $a = \infty$ ($T_{a,b} = T_b$). Moreover, the same argument proves $bP\{T_b \leq \delta_0 b/\mu\} = o(1)$ for all $\delta_0 < 1$.

Let $M^*$ denote universal constants. For $c_n \geq 1/n$, we split the integral for $\xi_j$ in (3.6) according to whether $H(x) \leq c_n$ so that

$$\xi_j = \Big(\int_{H \leq c_n} + \int_{H > c_n}\Big)\Big\{\log\Big(\frac{\psi_j}{\psi}\Big)d(jH_j) - \Big(\frac{jH_j}{H} - \frac{jW_j}{W}\Big)dH\Big\} = \xi_{n,j}^{(1)} + \xi_{n,j}^{(2)}.$$

Since $\psi_j = (jH_j)/(jW_j)$ with $jH_j$ and $jW_j$ being sums of iid measures, $\xi_{n,j+1}^{(1)} = \xi_{n,j}^{(1)}$ if $\int_{H \leq c_n} d\{(j+1)H_{j+1} - jH_j\} = 0$. Thus, by algebra and the Taylor expansion of $\log(\psi_{j+1}/\psi_j)$,

$$|\xi_{n,j+1}^{(1)} - \xi_{n,j}^{(1)}| \leq M^* \int_{H \leq c_n} \{1 + \log n - \log(H/2)\}d\{(j+1)H_{j+1} - jH_j\}$$

due to $d\{(j+1)W_{j+1} - jW_j\} \leq (1 \vee \Delta)d\{(j+1)H_{j+1} - jH_j\}$ and $\int (1/H_n)dH_n = \sum_{k=1}^{2n} 1/k$. For $\xi_{n,j}^{(2)}$, the four term Taylor expansion of $\log(\psi_j/\psi)$ yields

$$E \max_{n \leq j \leq 2n} \Big|\xi_{n,j}^{(2)} - U_{n,j}^{[1]} - (j-1)U_{n,j}^{[2]}\Big|^2 \leq M^*/(nc_n),$$

where $\{U_{n,j}^{[k]}, j \geq 1\}$ are $U$-statistics with kernels $U_{n,k}^{[k]}$ such that

$$E\Big|U_{n,1}^{[1]}\Big| \leq M^* \log\Big(\frac{2}{c_n}\Big), \quad \mathrm{Var}\Big(U_{n,1}^{[1]}\Big) \leq \frac{M^*}{nc_n}, \quad E\Big(U_{n,2}^{[2]}\Big)^2 \leq M^* \log\Big(\frac{2}{c_n}\Big),$$

with completely degenerate $U_{n,2}^{[2]}$. Choose $c_n$ satisfying $(\log n)c_n n^\alpha + 1/(nc_n) = o(1)$, we find that the above five inequalities imply the 1-regularity of $\xi_n$ for $p = 1$ and all $1/2 < \alpha < 1$.

**Proposition 3.1.** *The rank log-likelihood ratio $Z_n$ in (3.1) and (3.3) can be written as a perturbed random walk $Z_n = S_n + \xi_n$ with the random walk $S_n$ in (3.5), drift $\mu$ in (3.2), and the perturbation $\xi_n$ in (3.6). Moreover, for all $1/2 < \alpha < 1$ and $\epsilon > 0$, $E\max_{j \leq n^\alpha}|\xi_{n+j} - \xi_n| = o(1)$, $E(|\xi_n| - n^\alpha)^+ = o(1)$, $nP\{\max_{0 \leq j \leq n}|\xi_{n+j}| > \epsilon n^\alpha\} = o(1)$, and $\sum_n P\{|\xi_n| > \epsilon n^\alpha\} < \infty$.*

It follows from (3.5) that $S_1$ is non-lattice with $E|S_1|^p < \infty$ for all $p$. Since we have already proved $bP\{T_b \leq \delta_0 b/\mu\} = o(1)$, all conditions of Theorem 2.1 hold by Proposition 3.1. Thus, (3.7) holds and it remains to compute $E\xi_n$ up to $o(1)$. We sketch below the calculation in the case of $F = G$.

Assume $F = G$ are uniform in $(0,1)$. Let $u_1 \leq \cdots \leq u_{2n}$ be the ordered observations in the combined sample. Let $\varepsilon_k = 1$ if $u_k$ is from $G$ and $\varepsilon_k = 0$ otherwise. Set $y_k = \sum_{i=1}^k \varepsilon_i$. Since $y_k = nG_n(u_k)$ and $k = nH_n(u_k)$, $(\psi/\psi_n)(u_k) = 1 + \eta(2y_k/k - 1)$, so that by (3.6)

$$(3.8) \quad \xi_n = -\sum_{k=1}^{2n} \log\Big(1 + \eta(2y_k/k - 1)\Big) + \sum_{k=1}^{2n} \eta(2y_k - k)\log(u_{k+1}/u_k)$$

with $u_{2n+1} = 1$. Since $y_k \sim$ hypergeometric$(n, n, k)$, $Ey_k = k/2$. This and the independence of $\{y_k\}$ and $\{u_k\}$ imply

$$E\xi_n = (\eta^2/2)\log(2n) - C(\eta) + o(1).$$

Thus, the expansion in Theorem 3.1 holds for $A = 1$, in view of (3.7).

Since $(u_k - u_{k-1})/u_k$ are independent beta$(1, k-1)$ variables, (3.8) implies

$$\begin{aligned}
\xi_n &\approx \sum_{k=1}^{2n} \frac{\eta^2}{2}\left(\frac{2y_k}{k} - 1\right)^2 + \sum_{k=1}^{2n} \eta\left(\frac{2y_k}{k} - 1\right)\left\{(k+1)\frac{u_{k+1} - u_k}{u_{k+1}} - 1\right\} \\
&\approx \frac{\eta^2}{2}\int_{1/(2n)}^{1} \frac{B_0^2(t)}{t^2}dt + \eta \int_{1/(2n)}^{1} \frac{B_0(t)}{t}dW_1(t) \\
&\approx E\xi_n + N\Big(0, (\eta^4 + \eta^2)\log n\Big),
\end{aligned}$$
(3.9)

where $B_0(\cdot)$ is a Brownian bridge, $W_1(\cdot)$ is a Brownian motion independent of $B_0(\cdot)$, and $\xi_n \approx \xi'_n$ means $\xi_n$ has the same limiting distribution as $\xi'_n$ after centering by $E\xi_n$ and normalization by $\sqrt{\log n}$. Thus, condition (1.4) of the existing nonlinear renewal theorems does not hold.

## 4. Uniform integrability and an intermediate nonlinear renewal theorem

The results in this section are obtained by comparing $T_b$ in (1.2) with

$$\tau_b^* = \inf\{n \geq n_* : S_n + \zeta_{n_*} > b\}, \tag{4.1}$$

where $\zeta_n$ is as in (2.4) and $n_* = \lfloor b/\mu - \eta_* b^\alpha \rfloor$ for certain $\eta_* > \theta/\mu^{1+\alpha}$, $\eta_* < 1/\mu$ for $\alpha = 1$. We shall state all the results before proceeding to the proofs.

**Theorem 4.1.** *Suppose that $\{\xi_n, n \geq 1\}$ is $\rho$-regular with parameters $p \geq 1$ and $1/2 < \alpha \leq 1$. Suppose that the stopping time $T_b$ in (1.2) satisfies*

$$\lim_{b \to \infty} \Big(\rho(b)/b\Big)^{-p} P\{T_b \leq \delta_0 b/\mu\} = 0 \tag{4.2}$$

*for the $\delta_0$ in (2.2). If $E|X|^{(p+1)/\alpha} < \infty$, then*

$$\left\{\left|\frac{T_b - \tau_b^*}{\rho(b)}\right|^p : b \geq 1\right\} \quad \text{is uniformly integrable.} \tag{4.3}$$

For $\rho(x) = 1$ and $p = 1$, this uniform integrability result is a crucial component in our proof of the second order expansion of $ET_b$ in Theorem 2.1. It can be also used to derive expansions of the expectations of the renewal measure $U_b$ and the last exit time $N_b^*$, where

$$U_b = \sum_{n=1}^{\infty} I\{Z_n \leq b\}, \qquad N_b^* = 1 + \sup\{n \geq 1 : Z_n \leq b\}. \tag{4.4}$$

For general $\rho(\cdot)$, it yields an expansion of $ET_b$ up to $O(\rho(b))$. Since $\rho(b) = o(b)$ (e.g. $\rho(b)/b^\epsilon \to 0$ for all $\epsilon > 0$ with the rank SPRT by Lai [14], as shown in Section 3), such an expansion is typically sharper than direct extensions of the elementary renewal theorem but cruder than Theorem 2.1 as an extension of the standard renewal theorem. Thus, we call Theorem 4.2 (i) below an intermediate nonlinear renewal theorem.

**Theorem 4.2.** (i) *Suppose the conditions of Theorem 4.1 hold with $p = 1$. Then, for certain finite $M$,*

$$
\begin{aligned}
(4.5) \quad \mu E T_b &= b - E\zeta_{n_b} + O(1)\Big\{1 + E \max_{1 \leq j \leq Mb^\alpha} \big(|\xi_{n_*+j} - \zeta_{n_*}| \wedge n^\alpha\big)\Big\} \\
&= b - E\zeta_{n_b} + O(\rho(b)).
\end{aligned}
$$

(ii) *Suppose the conditions of Theorem 4.1 hold with $p = 2$. Then,*

$$(4.6) \quad \mathrm{Var}(T_b) = \frac{\sigma^2 b}{\mu^3} + O(\sqrt{b}\rho(b) + \rho^2(b)), \quad \text{as } b \to \infty.$$

**Remark 4.1.** The conclusions of Theorem 4.2 remain valid if we replace the $o(1)$ in (2.2) and (2.3) with $O(1)$ and (2.4) with the weaker

$$E \max_{1 \leq j \leq Mn^\alpha} \big(|\zeta_n - \xi_{n+j}|^p \wedge n^{\alpha p}\big) = O(\rho^p(b)).$$

**Theorem 4.3.** *Let $U_b$ and $N_b^*$ be as in (4.4). Suppose*

$$(4.7) \quad \sum_{k \geq n + Kn^\alpha} k^{p-1} P\Big\{\sup_{j \geq k} j^{-\alpha}(\xi_j + (j-n)\mu) \leq w_0\Big\} = o(\rho^p(n)),$$

*for some $K > 0$ and $w_0 > 0$ as in (2.3). Then, (4.3), (4.5) and (4.6) hold with $T_b$ replaced by either $U_b$ or $N_b^*$ under their respective conditions in Theorems 4.1 and 4.2.*

Proofs of of Theorems 4.2 and 4.3 are omitted since they follow from standard methods in nonlinear renewal theory, cf. [28] or [32]. We need three lemmas for the proof of Theorem 4.1.

**Lemma 4.1.** *Suppose $E|X|^{(p+1)/\alpha} < \infty$. If (2.2) and (4.2) hold for certain $\theta > 0$, then $P\{T_b \leq b/\mu - \eta_* b^\alpha\} = o(\rho(b)/b)^p$ for all $\eta_* > \theta/\mu^{1+\alpha}$.*

*Proof.* Let $w > 0$ satisfy $(1-w)\mu^{1+\alpha}\eta_* = \theta$ with the $\theta$ in (2.2). Let $b > 1$ and $n_* = \lfloor b/\mu - \eta_* b^\alpha \rfloor$. Let $\delta_0$ be as in (2.2) and (4.2). By (1.2)

$$
\begin{aligned}
(4.8) \quad &P\{\delta_0 b/\mu < T_b \leq b/\mu - \eta_* b^\alpha\} \\
&\leq P\{S_n + \xi_n > b \quad \text{for some} \quad \delta_0 b/\mu < n \leq n_*\} \\
&\leq P\Big\{\max_{1 \leq n \leq n_*} S_n - n\mu > w\mu\eta_* b^\alpha\Big\} \\
&\quad + P\Big\{\max_{\delta_0 n_* < n \leq n_*} (\xi_n + n\mu) > b - w\mu\eta_* b^\alpha\Big\}.
\end{aligned}
$$

Since $E|X|^{(p+1)/\alpha} < \infty$, by Theorem 1 of Chow and Lai [5]

$$(4.9) \quad P\Big\{\max_{1 \leq n \leq b} |S_n - n\mu| > \lambda b^\alpha\Big\} = o(b^{-p}), \quad \forall \lambda > 0.$$

Thus, the first term on the right-hand side of (4.8) is of the order $o(b^{-p})$ since $n_* = O(b)$. Since $b - w\mu\eta_* b^\alpha - n_*\mu \geq (1-w)\mu\eta_* b^\alpha \geq (1-w)\mu^{1+\alpha}\eta_* n_*^\alpha = \theta n_*^\alpha$ and $\rho(n_*)$ is of the same order as $\rho(b)$ by (2.1), it follows from (2.2) that the second term on the right-hand side of (4.8) is bounded by

$$P\Big\{\max_{\delta_0 n_* < n \leq n_*} \xi_n > \theta n_*^\alpha\Big\} = o(\rho(n_*)/n_*)^p = o\big(\rho(b)/b\big)^p.$$

Hence, (4.8) is $o\big(\rho(b)/b\big)^p$ and the conclusion follows from (4.2). □

**Lemma 4.2.** *Suppose $E|X|^{(p+1)/\alpha} < \infty$. If (2.3) holds for certain $K > 0$, then*

$$\lim_{b\to\infty} \frac{1}{\rho^p(b)} \sum_{k \geq b/\mu + \eta^* b^\alpha} k^{p-1} P\{T_b > k\} = 0, \quad \eta^* = K/\mu^\alpha.$$

*Proof.* Let $n^* = \lceil b/\mu + \eta^* b^\alpha \rceil$. By (1.2)

$$(4.10) \quad \sum_{k=n^*}^{\infty} k^{p-1} P\{T_b > k\} \leq \sum_{k=n^*}^{\infty} k^{p-1} P\{S_k + \xi_k \leq b\}$$

$$\leq \sum_{k=n^*}^{\infty} k^{p-1} P\{S_k \leq k\mu - w_0 k^\alpha/2\}$$

$$+ \sum_{k=n^*}^{\infty} k^{p-1} P\{\xi_k \leq b - k\mu + w_0 k^\alpha/2\}.$$

Since $E|X|^{(p+1)/\alpha} < \infty$, by Theorem 1 of Chow and Lai [5]

$$(4.11) \quad \sum_{n=1}^{\infty} P\{|S_n - n\mu| > \lambda n^\alpha\} < \infty, \quad \forall \lambda > 0.$$

Thus, since $w_0 > 0$, the first term on the right-hand side of (4.10) is $o(1)$ as $b \to \infty$. Let $n_b = \lfloor b/\mu \rfloor$ as in Theorem 2.1. Since $n_b\mu \leq b < (n_b + 1)\mu$ and $n^* \geq n_b + \eta^*(n_b\mu)^\alpha = n_b + Kn_b^\alpha$, in view of (2.3) the second term is bounded by

$$\sum_{k \geq n_b + Kn_b^\alpha} k^{p-1} P\{\xi_k \leq (n_b+1)\mu - k\mu + w_0 k^\alpha/2\} = o(\rho^p(n_b)).$$

The conclusion follows since $\rho(n_b)$ and $\rho(b)$ are of the same order. □

**Lemma 4.3.** *Suppose $E|X|^{(p+1)/\alpha} < \infty$. Let $\tau_b^*$ be as in (4.1) with the $\eta_*$ in Lemma 4.1. Then,*

$$(4.12) \quad P\{b/\mu - \eta_* b^\alpha < \tau_b^* < b/\mu + \eta_\tau^* b^\alpha\} \geq 1 + o(b^{-p})$$

*for all $\eta_\tau^* > \theta^*/\mu^{1+\alpha}$. Moreover, for such $\eta_\tau^*$*

$$(4.13) \quad \lim_{b\to\infty} \sum_{k \geq b/\mu + \eta_\tau^* b^\alpha} k^{p-1} P\{\tau_b^* > k\} = 0.$$

*Proof.* Since $b - n_*\mu - \theta n_*^\alpha \geq \mu\eta_* b^\alpha - \theta n_*^\alpha \geq (\eta_*\mu^{1+\alpha} - \theta)n_*^\alpha > 0$, it follows from (4.9) that

$$P\{\tau_b^* = n_*\} \leq P\{S_{n_*} + \theta n_*^\alpha > b\}$$
$$= P\{S_{n_*} - n_*\mu > b - n_*\mu - \theta n_*^\alpha\} = o(b^{-p}).$$

Let $n_\tau^* = \lceil b/\mu + \eta_\tau^* b^\alpha \rceil$. Since $b - n_\tau^*\mu + \theta^* n_*^\alpha \leq -\mu\eta_\tau^* b^\alpha + \theta^* n_*^\alpha \leq -(\eta_\tau^*\mu^{1+\alpha} - \theta^*)n_*^\alpha < 0$, we have

$$P\{\tau_b^* \geq n_\tau^*\} \leq P\{S_{n_\tau^*} - \theta^* n_*^\alpha \leq b\}$$
$$= P\{S_{n_\tau^*} - n_\tau^*\mu \leq b - n_\tau^*\mu + \theta^* n_*^\alpha\} = o(b^{-p}).$$

The above calculations prove (4.12). The proof of (4.13) is simpler than that of Lemma 4.2 and omitted. □

*Proof of Theorem* 4.1. Let $\eta_*, \eta^* = K/\mu^\alpha$ and $\eta^* > \eta_\tau^* > \theta^*/\mu^{1+\alpha}$ be as in Lemmas 4.1, 4.2 and 4.3 respectively. This is possible since $\theta^* < K\mu$ with the $\rho$-regularity conditions. Set $n_* = \lfloor b/\mu - \eta_* b^\alpha \rfloor$, $n^* = \lceil b/\mu + \eta^* b^\alpha \rceil$ and $n_\tau^* = \lceil b/\mu + \eta_\tau^* b^\alpha \rceil$.

*Step 1:* We first prove that for all integers $k_* \geq 1$

$$\sum_{k=k_*}^\infty k^{p-1} P\{T_b - \tau_b^* > k\}$$

(4.14)
$$\leq 2 \sum_{k \geq n^*} k^{p-1} P\{T_b > k\} + \sum_{n^*-n_* \leq k < n^*} k^{p-1} P\{\tau_b^* \leq n_*\}$$
$$+ \sum_{k_* \leq k < n^*-n_*} k^{p-1} P\{A_\tau\} + \sum_{k_* \leq k < n^*-n_*} k^{p-1} P\{S_k \leq \mu k/2\}$$
$$+ \sum_{k_* \leq k < n^*-n_*} k^{p-1} P\Big\{\max_{1 \leq j \leq 2+M_\tau b^\alpha}(\zeta_{n_*} - \xi_{j+n_*}) > \mu k/2\Big\}$$

with $A_\tau = \{\tau_b^* \notin (n_*, n_\tau^*)\}$ and $M_\tau = \eta_\tau^* + \eta^* + 2\eta_*$, and

$$\sum_{k=k_*}^\infty k^{p-1} P\{\tau_b^* - T_b > k\}$$

(4.15)
$$\leq 2 \sum_{k \geq n^*} k^{p-1} P\{\tau_b^* > k\} + \sum_{n^*-n_* \leq k < n^*} k^{p-1} P\{T_b \leq n_*\}$$
$$+ \sum_{k_* \leq k < n^*-n_*} k^{p-1} P\{\widetilde{A}_T\} + \sum_{k_* \leq k < n^*-n_*} k^{p-1} P\{S_k \leq \mu k/2\}$$
$$+ \sum_{k_* \leq k < n^*-n_*} k^{p-1} P\Big\{\max_{1 \leq j \leq 2+M_T b^\alpha}(\xi_{j+n_*} - \zeta_{n_*}) > \mu k/2\Big\}$$

with $\widetilde{A}_T = \{T_b \notin (n_*, n^*), \tau_b^* > T_b\}$ and $M_T = \eta^* + \eta_*$.

For the proof of (4.14), we divide $[k_*, \infty)$ into three intervals at $n^* - n_*$ and $n^*$. For the first interval $k \in [k_*, n^* - n_*)$, we have $n_* < \tau_b^* + k \leq n_\tau^* + n^* - n_* - 2 < n_* + 2 + M_\tau b^\alpha$ in $A_\tau^c$ and

$$\{T_b - \tau_b^* > k\} \subseteq \{S_{\tau_b^*+k} + \xi_{\tau_b^*+k} \leq b < S_{\tau_b^*} + \zeta_{n_*}\}$$
$$\subseteq \{S_{\tau_b^*+k} - S_{\tau_b^*} \leq \mu k/2\} \cup \{\xi_{\tau_b^*+k} - \zeta_{n_*} < -\mu k/2\}.$$

For the second interval $k \in [n^* - n_*, n^*)$, we have

$$\sum_{k \geq n^*-n_*} k^{p-1} P\{T_b > k + \tau_b^*, \tau_b^* > n_*\} \leq \sum_{k \geq n^*-n_*} k^{p-1} P\{T_b > k + n_*\}$$
$$\leq \sum_{j \geq n^*} j^{p-1} P\{T_b \geq j\}.$$

The proof of (4.15) is nearly identical, with $n_* < T_b < n^* \leq n_* + (\eta^* + \eta_*)b^\alpha + 2$ in $\widetilde{A}_T^c$ and

$$\{\tau_b^* - T_b > k\} \subseteq \{S_{T_b+k} + \zeta_{n_*} \leq b < S_{T_b} + \xi_{T_b}\}$$
$$\subseteq \{S_{T_b+k} - S_{T_b} \leq \mu k/2\} \cup \{\zeta_{n_*} - \xi_{T_b} < -\mu k/2\}.$$

*Step 2:* Prove the uniform integrability of $\{(T_b - \tau_b^*)^+/\rho(b)\}^p$. For $k_* = \lfloor C\rho(b) \rfloor$, we have

$$J(b, C) = E\left(\frac{T_b - \tau_b^*}{\rho(b)} - C\right)^p I\left\{\frac{T_b - \tau_b^*}{\rho(b)} > C\right\}$$

$$\leq \frac{C_p}{\rho^p(b)} \sum_{k=k_*}^{\infty} k^{p-1} P\{T_b - \tau_b^* > k\},$$

where $C_p$ is a universal constant. For the current step, it suffices to show $J(b,C) \to 0$ as $b \to \infty$ and then $C \to \infty$. Since $n^* = O(b)$ and $n^* - n_* = O(b^\alpha)$, it follows from (4.14), Lemmas 4.2 and 4.3, and (4.11) that

$$J(b,C) = o(1) + \frac{C_p}{\rho^p(b)} \sum_{k_* \leq k < n^* - n_*} k^{p-1} P\Big\{ \max_{1 \leq j \leq 2 + M_\tau b^\alpha} (\zeta_{n_*} - \xi_{j+n_*}) > \mu k/2 \Big\}.$$

Since $b = O(n_*)$, we may choose $M$ in (2.4) satisfying $Mn_*^\alpha > 2 + M_\tau b^\alpha$ for all $b > 1$. Thus, by the uniform integrability in (2.4), $J(b,C) = o(1)$ as $b \to \infty$ and then $C \to \infty$.

*Step 3:* Prove the uniform integrability of $\{(T_b - \tau_b^*)^-/\rho(b)\}^p$. This step is nearly identical to Step 2, with (4.15) instead of (4.14). In fact, since $M_T < M_\tau$, the same $M$ in (2.4) works. Although Lemma 4.1 does not provide $(b^\alpha/\rho(b))^p P(T_b \geq n^*) = o(1)$, we have $P(\widetilde{A}_T) \leq P(T_n \leq n_*) + P(\tau_b^* \geq n^* \geq n_\tau^*)$, so that Lemmas 4.1 and 4.3 can be used to control the third sum on the right-hand side of (4.15). This completes the proof of Step 3 and thus the entire theorem. □